\documentclass[11pt]{article}
\usepackage{graphicx}
\usepackage[compress]{cite}
\usepackage{amsmath}
\usepackage{amsfonts}
\usepackage[papersize={8.5in,11in},margin=1in]{geometry}

\newcommand{\EQ}[1]{(\ref{equation:#1})}
\DeclareMathOperator*{\argmin}{arg\,min}

\begin{document}
\begin{center}
\Large \bf{Parameter estimation by implicit sampling} 
\end{center}
\vspace{1mm}

\begin{center}
  Matthias Morzfeld$^{1,2,*}$, Xuemin Tu$^{3}$,  Jon Wilkening$^{1,2}$ and Alexandre J. Chorin$^{1,2}$

\vspace{3mm}
$^1$Department of Mathematics, University of California, Berkeley, CA 94720,~USA.\\
$^2$Lawrence Berkeley National Laboratory, Berkeley, CA 94720,~USA.\\
$^3$Department of Mathematics, University of Kansas, Lawrence, KS 66045,~USA.\\
\vspace{1mm}
\let\thefootnote\relax\footnote{$^*$ Corresponding author. Tel: +1~510~486~6335. Email address: mmo@math.lbl.gov. 
Lawrence Berkeley National Laboratiry, 1 Cyclotron Road, Berkeley, California 94720,~USA.}

\emph{Abstract}
\end{center}
Implicit sampling is a weighted sampling method that is used in data assimilation, 
where one sequentially updates estimates of the state of a stochastic model
based on a stream of noisy or incomplete data.
Here we describe how to use implicit sampling in parameter estimation problems,
where the goal is to find parameters of a numerical model, e.g.~a 
partial differential equation (PDE), 
such that the output of the numerical model is compatible with (noisy) data. 
We use the Bayesian approach to parameter estimation, in which a posterior probability density 
describes the probability of the parameter  conditioned on data
and compute an empirical estimate of this posterior with implicit sampling.
Our approach generates independent samples,
so that some of the practical difficulties one encounters with Markov Chain Monte Carlo methods,
e.g.~burn-in time or correlations among dependent samples, are avoided.
We describe a new implementation of implicit sampling for parameter estimation problems
that makes use of multiple grids (coarse to fine) and BFGS optimization coupled to adjoint equations
for the required gradient calculations.
The implementation is ``dimension independent'',
in the sense that a well-defined finite dimensional subspace is sampled
as the mesh used for discretization of the PDE is refined.
We illustrate the algorithm with an example
where we estimate a diffusion coefficient in an elliptic equation
from sparse and noisy pressure measurements.
In the example, dimension\slash mesh-independence is achieved via Karhunen-Lo\`{e}ve expansions.

\section{Introduction}
We take the Bayesian approach to 
parameter estimation and compute the probability density function (pdf)
$p(\theta|z)$, where $\theta$ is a set of parameters (an $m$-dimensional vector)
and $z$ are data (a $k$-dimensional vector, see, e.g.~\cite{StuartInverse}).
We assume a prior pdf $p(\theta)$ for the parameters,
which describes what one knows about the parameters before 
collecting the data.
For example, one may know a priori that a parameter is positive.
We further assume a likelihood  $p(z|\theta)$, which describes
how the parameters are connected with the data.
Bayes' rule combines the prior and likelihood to find 
$p(\theta|z)\propto p(\theta)p(z|\theta)$ as a posterior density. 

If the prior and likelihood are Gaussian, 
then the posterior is also Gaussian,
and it is sufficient to compute the
mean and covariance of $\theta|z$
(because the mean and covariance define the Gaussian).
Moreover, the mean and covariance
are the minimizer and the inverse of the Hessian of
the negative logarithm of the posterior.
If the posterior is not Gaussian, 
e.g.~because the numerical model is nonlinear,
then one can compute the posterior mode,
often called the maximum a posteriori (MAP) point, 
by minimizing the negative logarithm of the posterior,
and use the MAP point
as an approximation of the parameter $\theta$.
The inverse of the Hessian of the negative
logarithm of the posterior can be used to measure  
the uncertainty of the approximation.
This method is sometimes called linearization about the MAP point (LMAP)
or the Laplace approximation~\cite{Iglesias2,Oliver2011,Dean2007,Bui}.

One can also use Markov Chain Monte Carlo (MCMC) 
and represent the posterior by a collection of samples,
see, e.g.~\cite{Efendiev,Martin,Petra,Stuart2}.
The samples form an empirical estimate of the posterior,
and statistics, e.g.~the mean or mode, 
can be computed from this empirical estimate
by averaging over the samples.
Under mild assumptions, the moments one computes from the samples
converge to the moments of the posterior
(as the number of samples goes to infinity).
In practice, one often encounters difficulties with MCMC.
For example, the samples may have a distribution
which converges slowly to the posterior, 
or there could be strong correlations among the samples.
In these cases, it is difficult to determine
how many samples are ``enough samples'' for an accurate empirical estimate.

We propose to use importance sampling to avoid some of these issues,
in particular the estimation of burn-in or correlation times.
The idea in importance sampling is to generate independent 
samples from a density that one knows how to sample, the importance function,
rather than from the one one wants to sample.
Weights are attached to each sample to
account for the imperfection of the importance function.
Under mild assumptions, the weighted samples also form an 
empirical estimate of the posterior pdf~\cite{ChorinHald}.
The efficiency and applicability of an importance sampling scheme
depends on the choice of the importance function. 
Specifically, a poorly chosen importance function can be (nearly) singular
with respect to the posterior,
in which case most of the samples one generates 
are unlikely with respect to the posterior,
so that
the number of samples required 
becomes large and importance sampling, therefore, becomes impractical
\cite{Blb,Blb2,Sny}.

We show how to use implicit sampling 
for constructing importance functions that are
large where the posterior is large,
so that a manageable number of samples forms an accurate empirical estimate.
Implicit sampling has been studied before in the context of
online-filtering\slash data assimilation, i.e.~state estimation of a stochastic
model in \cite{amc,Morzfeld2012,cmt,CT1,Morzfeld2011,cmt2},
and for parameter estimation in stochastic models in \cite{Brad}.
Here we describe how to use implicit sampling
for Bayesian parameter estimation.
In principle, using implicit sampling for parameter estimation is straightforward
(since it is a technique for sampling arbitrary, finite-dimensional probability densities~\cite{amc}),
however its implementation in the context of parameter estimation requires attention.
We discuss how to sample the posterior of parameter estimation
with implicit sampling (in general), 
as well as a specific implementation that
is suitable for parameter estimation,
where the underlying model
that connects the parameters with the data is
a partial differential equation (PDE).
We show that the sampling algorithm is 
independent of the mesh one uses for discretization of the PDE.
The idea of mesh-independence is also discussed in,e.g.,~\cite{Bui,Martin,Petra}.
Mesh-independence means that the sampling algorithm
``converges'' as the mesh is refined in the sense that
the same finite dimensional subspace is sampled.
We further show how to use multiple grids
and adjoints during the required optimization
and during sampling, and discuss approximations
of a Hessian to reduce the computational costs.
The optimization in implicit sampling represents the link 
between implicit sampling and LMAP, as well as recently 
developed stochastic Newton MCMC methods~\cite{Bui,Martin,Petra}.
The optimization and Hessian codes used in these codes can be used
for implicit sampling, and the weighing of the samples 
can describe non-Gaussian characteristics of the posterior
(which are missed by LMAP).

We illustrate the efficiency of our implicit sampling 
algorithm with numerical experiments
in which we estimate the diffusion coefficient in 
an elliptic equation using sparse and noisy data. 
This problem is a common test problem for MCMC algorithms
and has important applications in reservoir simulation\slash management and in pollution modeling \cite{Bear,Dean2007}.
Moreover, the conditions for the existence of a posterior measure and its continuity are well understood \cite{Stuart2}. 
Earlier work on this problem includes \cite{Efendiev}, where Metropolis-Hastings MC sampling is used, 
\cite{Iglesias} where an ensemble Kalman filter is used, 
and~\cite{MarzoukOptimalMaps}, which uses optimal maps and is further discussed below. 

The remainder of this paper is organized as follows. 
In section 2 we explain how to use implicit sampling 
for parameter estimation and discuss an efficient implementation.
A numerical example is provided in section 3.
The numerical example involves an elliptic PDE
so that the dimension of the parameter we estimate is infinite.
We discuss its finite dimensional approximation 
and achieve mesh-independence via KL expansions. 
Conclusions are offered in section 4.

\section{Implicit sampling for parameter estimation}
We wish to estimate an $m$-dimensional parameter vector $\theta$ from data 
which are obtained as follows. 
One measures a function of the parameters $h(\theta)$, 
where $h$ is a given $k$-dimensional function; 
the measurements are noisy, so that the data $z$ satisfy the relation:
\begin{equation}
\label{equation:obs}
z=h(\theta)+r,
\end{equation}
where $r$ is a random variable with a known distribution 
and the function $h$ maps the parameters onto the data. 
Often, the function $h$ involves solving a PDE.
In a Bayesian approach, 
one obtains the pdf $p(\theta |z)$ of the conditional random 
variable $\theta | z$ by Bayes' rule:
\begin{equation}
	p(\theta|z)\propto p(\theta)p(z|\theta),
\label{bbe}
\end{equation}
where the likelihood $p(z|\theta)$ can be read off (\ref{equation:obs}) 
and the prior $p(\theta)$ is assumed to be known. 

The goal is to compute this posterior.
This can be done with importance sampling as follows~\cite{ChorinHald,Kalos}.
One can represent the posterior by $M$ weighted samples. 
The samples~$\theta_j$, $j=1,\dots,M$ are obtained 
from an importance function $\pi(\theta)$ 
(which is chosen such that it is easy to sample from),
and the $j$th sample is assigned the weight
\begin{equation*}
	w_j\propto \frac{p(\theta_j)p(z|\theta_j)}{\pi(\theta_j)}.
\end{equation*}
The location of a sample corresponds to a set of possible parameter values 
and the weight describes how likely this set is in view of the posterior. 
The weighted samples $\left\{ \theta_j,w_j \right\}$ form an empirical estimate of $p(\theta|z)$, 
so that for a smooth function $u$, the sum 
\begin{equation*}
	E_M(u) = \sum_{j=0}^{M}u(\theta_j)\hat{w}_j,
\end{equation*}
where $\hat{w}_j=w_j/\sum_{j=0}^{M}w_j$, 
converges almost surely to the expected value of $u$ 
with respect to $p(\theta|z)$ as $M \rightarrow \infty$, 
provided that the support of $\pi$ includes the support of~$p(\theta|z)$ \cite{Kalos,ChorinHald}. 

The importance function must be chosen carefully
or else sampling is inefficient \cite{Blb,Blb2,Sny,CM2013}.
We construct the importance function via implicit sampling, 
so that the importance function is large where the posterior pdf is large.
This is done by computing the maximizer of $p(\theta|z)$,
i.e.~the MAP point. 
If the prior and likelihood are exponential functions (as they often are in applications), 
it is convenient to find the MAP point by minimizing the function 
\begin{equation}
\label{equation:FDef}
	F(\theta) = -\log\left( p(\theta)p(z|\theta)\right).
\end{equation}

After the minimization, 
one finds samples in the neighborhood of the MAP point, $\mu = \argmin F$, as follows.
Choose a reference variable $\xi$ with pdf $g(\xi)$ and let $G(\xi)=-\log(g(\xi))$ and $\gamma=\min_\xi G$.
For each member of a sequence of samples of $\xi$ solve the equation
\begin{equation}
\label{gen}
F(\theta)-\phi=G(\xi)-\gamma,
\end{equation}
to obtain a sequence of samples $\theta$, where $\phi$ is the minimum of $F$. 
The sampling weight is 
\begin{equation}
w_j\propto J(\theta_j),
\end{equation}
where $J$ is the Jacobian of the one-to-one and onto map $\theta \rightarrow \xi$. 
There are many ways to choose this map since (\ref{gen}) is underdetermined \cite{cmt,CT1,Morzfeld2011}; 
we describe two choices below.  
The sequence of samples we obtain by solving (\ref{gen}) is in the neighborhood of the minimizer $\mu$ since, 
by construction, equation~(\ref{gen}) maps a likely $\xi$ to a likely $\theta$: 
the right hand side of (\ref{gen}) is small with a high probability since 
$\xi$ is likely to be close to the mode (the minimizer of $G$); 
thus the right hand side is also likely to be small and, therefore, 
the sample is in the neighborhood of the MAP point $\mu$. 

An interesting construction, related to implicit sampling,
has been proposed in \cite{MarzoukOptimalMaps,Tabak}. 
Suppose one wants to generate samples with the pdf $p(\theta|z)$, and have 
$\theta$ be a function of a reference variable $\xi$ with pdf $g$, as above.
If the samples are all to have equal weights, one must have, in the notations above,
\begin{equation*}
	p(\theta|z)=g(\xi)/J(\xi),
\end{equation*}
where, as above, $J$ is the Jacobian of a map $\theta \rightarrow \xi$.
Taking logs, one finds
\begin{equation}
\label{Marzouk}
F(\theta)+\log\beta=G(\xi)-\log\left(J(\xi)\right),
\end{equation}
where $\beta=\int p(z|\theta)p(\theta)d\theta$ is the 
proportionality constant that has been elided in~(\ref{bbe}). 
If one can find a one-to-one mapping from $\xi$ to $\theta$ that satisfies this equation, 
one obtains an optimal sampling strategy, 
where the pdf of the samples matches exactly the posterior pdf. 
In \cite{MarzoukOptimalMaps}, this map is found globally by choosing $g=p(\theta)$ (the prior), 
rather than sample-by-sample as in implicit sampling. 
The main differences between the implicit sampling equation (\ref{gen}) and equation~(\ref{Marzouk}) 
are the presence of the Jacobian $J$ and of the normalizing constant $\beta$ in the latter; 
$J$ has shifted from being a weight to being a term in the equation that picks the samples, 
and the optimization that finds the probability mass has shifted to the computation of the map.

If the reference variable is Gaussian and the problem is linear, 
equation~(\ref{Marzouk}) can be solved by a linear map with a constant Jacobian
and this map also solves (\ref{gen}), 
so that one recovers implicit sampling. 
In particular, in a linear Gaussian problem, 
the local (sample-by-sample) map (\ref{gen}) of implicit sampling 
also solves the global equation (\ref{Marzouk}), which, for the linear problem, 
is a change of variables from one Gaussian to another. 
If the problem is not linear, 
the task of finding a global map that satisfies (\ref{Marzouk}) 
is difficult (see also \cite{Doucet2000,liuchen1995,Zaritski,Tabak}). 
The determination of optimal maps in~\cite{MarzoukOptimalMaps}, 
based on nonlinear transport theory, 
is elegant but can be computationally intensive, 
and requires approximations that reintroduce non-uniform weights. 
Using (simplified) optimal maps and re-weighing the samples
from approximate maps is discussed in~\cite{Tabak}.
In \cite{Oliver14}, further optimal transport maps from
prior to posterior are discussed. 
These maps are exact in linear Gaussian problems,
however in general they are approximate, due to the neglect of a Jacobian,
when the problem is nonlinear.

\subsection{Implementation of implicit sampling for parameter estimation}
Above we assume that the parameter $\theta$ is finite-dimensional.
However, if $h$ involves a PDE, then the parameter may be infinite-dimensional.
In this case, one can discretize the parameter,
e.g.~using Karhunen-Lo\`{e}ve expansions (see below).
The theory of implicit sampling then immediately applies.
A related idea is dimension-\slash mesh-independent MCMC, see,e.g., 
\cite{StuartInverse,Stuart2,Bui,Petra,Martin},
where MCMC sampling algorithms operate efficiently,
independently of the mesh that is used to discretize the PDE.
In particular, the algorithms sample the same finite dimensional sub-space
as the mesh-size is refined.
Below we present an implementation of implicit sampling
that is mesh independent.

\subsubsection{Optimization and multiple-grids}
The first step in implicit sampling is to find the MAP point by
minimizing $F$ in (\ref{equation:FDef}).
Upon discretization of the PDE, 
this can be done numerically by Newton, quasi-Newton, or Gauss-Newton methods 
(see, e.g.~\cite{Nocedal}).
The minimization requires derivatives of the function $F$, 
and these derivatives may not be easy to compute.
When the function $h$ in~(\ref{equation:obs}) involves solving a PDE,
then adjoints are efficient for computing the gradient of $F$.
The reason is that the complexity of solving the adjoint equation is similar to that of 
solving the original ``forward'' model.
Thus, the gradient can be computed at the cost of (roughly) two forward solutions of (\ref{equation:obs}). 
The adjoint based gradient calculations can be used in connection with a quasi-Newton method, e.g.~BFGS,
or with Gauss-Newton methods.
We illustrate how to use the adjoint method for BFGS optimization in the example below.

One can make use of multiple grids during the optimization. 
This idea first appeared in the context of online state estimation in \cite{amc}, 
and is  similar to a multi-grid finite difference method \cite{Fedorenko} 
and multi-grid Monte Carlo~\cite{Goodman}. 
The idea is as follows. 
First, initialize the parameters and pick a coarse grid. 
Then perform the minimization on the coarse grid and use the minimizer 
to initialize a minimization on a finer grid. 
The minimization on the finer grid should require only a few steps, 
since the initial guess is informed by the computations on the coarser grid,
so that the number of fine-grid forward and adjoint solves is small. 
This procedure can be generalized to use more than two grids (see the example below).

For the minimization, we distinguish between two scenarios.
First, suppose that the physical parameter can have a mesh-independent representation
even if it is defined on a grid.
This happens, for example, if the parameter is represented using
Karhunen-Lo\`{e}ve expansions, where the expansions are later evaluated on the grid.
In this case, no interpolation of the parameter is required in the multiple-grid approach.
On the other hand, if the parameter is defined on the grid,
then the solution on the coarse mesh must be interpolated onto the fine grid,
as is typical in classical multi-grid algorithms (e.g.~for linear systems). 
We illustrate the multiple-grid approach in the context 
of the numerical example in section~3.

\subsubsection{Solving the implicit equations with  linear maps}
Once the optimization problem is done, 
one needs to solve the random algebraic equations (\ref{gen}) to generate samples. 
There are many ways to solve~(\ref{gen}) 
because it is an underdetermined equation (one equation in $m$ variables).
We describe and implement two strategies for solving~(\ref{gen}).
For a Gaussian reference variable $\xi$ with mean $0$ and covariance matrix $H^{-1}$,
where $H$ is the Hessian of the function $F$ at the minimum,
the equation becomes
\begin{equation}
\label{equation:SamplingEquation}
	F(\theta)-\phi =\frac{1}{2}\xi^TH\xi.
\end{equation}
In implicit sampling with linear maps
one approximates $F$ 
by its Taylor expansion to second order
\begin{equation}
	F_0(\theta) =\phi + \frac{1}{2}(\theta-\mu)^TH(\theta-\mu),
	\nonumber
\end{equation}
where $\mu=\mbox{arg min}\, F$ is the minimizer of $F$ (the MAP point) and 
$H$ is the Hessian at the minimum.
One can then solve the quadratic equation
\begin{equation}
	\label{eq:QuadApprox}
	F_0(\theta) -\phi = \frac{1}{2}\xi^TH\xi,
\end{equation}
instead of (\ref{equation:SamplingEquation}),
using
\begin{equation}
\label{eq:LM}
	\theta = \mu +  \xi.
\end{equation}
The bias created by solving the quadratic equation (\ref{eq:QuadApprox}) 
instead of (\ref{equation:SamplingEquation}) 
can be removed by the weights \cite{cmt,amc}
\begin{equation}
\label{eq:weightsLinear}
	w\propto \exp\left(F_0(\theta)-F(\theta)\right).
\end{equation}
Note that the algorithm is mesh-independent
in the sense of~\cite{Martin,Bui,Petra}
due to the use of Hessian of~$F$.
Specifically,
the eigenvectors associated with non-zero eigenvalues
of the discrete Hessian span the same
stochastic sub-space as the mesh is refined.

\subsubsection{Solving the implicit equations with random maps}
A second strategy for solving (\ref{gen}) is to use random maps~\cite{Morzfeld2011}. 
The idea is to solve (\ref{equation:SamplingEquation}) 
in a random direction, $\xi$, where $\xi\sim\mathcal{N}(0,H^{-1})$ as before:
\begin{equation}
\label{equation:Ansatz}
	\theta = \mu +\lambda(\xi)\,\xi.
\end{equation}
We look for a solution of (\ref{equation:SamplingEquation}) 
in a $\xi$-direction by substituting~(\ref{equation:Ansatz}) into~(\ref{equation:SamplingEquation}), 
and solving the resulting equation for the scalar~$\lambda(\xi)$ with Newton's method. 
A formula for the Jacobian of the random map defined by~(\ref{equation:SamplingEquation}) 
and~(\ref{equation:Ansatz}) was derived in \cite{Morzfeld2011,Goodman},
\begin{equation}
\label{eq:Jacobian}
	w \propto \left\vert J(\xi)\right\vert = \left\vert\lambda^{m-1} \;\frac{\xi^TH\xi}{\nabla_\theta F \cdot\xi}\right\vert
\end{equation}
where $m$ is the number of non-zero eigenvalues of $H$,
making it easy to evaluate the weights of the samples if the gradient of $F$ is easy to compute,
e.g.~using the adjoint method (see below).
Note that the random map algorithm is affine invariant and, therefore,
capable of sampling within flat and narrow valleys of $F$.
It is also mesh-independent in the sense of~~\cite{Martin,Bui,Petra},
for the same reasons as the linear map method above.

\section{Application to subsurface flow}
We illustrate the applicability of our implicit sampling method
by a numerical example from subsurface flow,
where we estimate subsurface structures from pressure measurements of flow through a porous medium.
This is a common test problem for MCMC and has applications in 
reservoir simulation\slash management (see e.g.~\cite{Dean2007}) 
and pollution modeling (see e.g.~\cite{Bear}). 

We consider Darcy's law
\begin{equation*}
	\kappa \nabla p = - \nu u,
\end{equation*}
where $\nabla p$ is the pressure gradient across the porous medium, 
$\nu$ is the viscosity and $u$ is the average flow velocity; 
$\kappa$ is the permeability and describes the subsurface structures we are interested in. 
Assuming, for simplicity, that the viscosity and density are constant, 
we obtain, from conservation of mass, the elliptic problem
\begin{equation}\label{equation:upeqn}
-\nabla\cdot \left(\kappa\nabla p\right) = g,
\end{equation}
on a domain $\Omega$, 
where the source term $g$ represents externally prescribed inward or outward flow rates. 
For example, if a hole were drilled and a constant inflow were applied through this hole, 
$g$ would be a delta function with support at the hole. 
Here we choose $g=200\pi^2 \sin(\pi x)\sin(\pi y)$.
Equation (\ref{equation:upeqn}) is supplemented with Dirichlet boundary conditions.

The uncertain quantity in this problem is the permeability, 
i.e.~$\kappa$ is a random variable, 
whose realizations we assume to be smooth enough 
so that for each one a solution of \EQ{upeqn} uniquely exists. 
We would like to update our knowledge about $\kappa$ on the basis of 
noisy measurements of the pressure at $k$ locations 
within the domain~$\Omega$ so that \EQ{obs} becomes
\begin{equation}\label{equation:obseqn}
z = h(p(\kappa),x,y) + r.
\end{equation}
Computation requires a discretization of the forward problem \EQ{upeqn} 
as well as a characterization of the uncertainty in the permeability before data are collected, 
i.e.~a prior for $\kappa$. 
We describe our choices for the discretization and prior below.

\subsection{Discretization of the forward problem}
In the numerical experiments below we consider a 2D-problem
and choose the domain $\Omega$ to be the square $[0,1]\times [0,1]$,
and discretize \EQ{upeqn} with a piecewise linear
finite element method on a uniform $(N+1)\times (N+1)$ mesh of triangular elements
with $2(N+1)^2$ triangles~\cite{Braess}.
Solving the (discretized) PDE thus
amounts to solving the linear system
\begin{equation}\label{equation:DEQ}
	AP=G,
\end{equation}
where $A$ is a $N^2\times N^2$ matrix,
and where $P$ and $G$
are $N^2$ vectors;
$P$ is the pressure and $G$ contains the discretized right
hand side of the equation \EQ{upeqn}.
For a given permeability~$\kappa$, the matrix $A$ is
symmetric positive definite (SPD) and 
we use the balancing domain decomposition by constraints method
\cite{BDDC} 
to solve (\ref{equation:DEQ}), 
i.e.~we first decompose the computational domain into smaller subdomains 
and then solve a subdomain interface problem. 
For details of the linear solvers, see~\cite{BDDC}. 

In the numerical experiments below,
a $64\times 64$ grid is our finest grid,
and the data, i.e.~the pressure measurements,
are arranged such that they align with grid points 
of our finest grids, as well as with the coarse grids
we use in the multiple-grid approach.
Specifically, the data equation~\EQ{obseqn} becomes 
\begin{equation*}
	z = MP + r,
\end{equation*}
where $M$ is a $k\times N^2$ matrix
that defines at which locations on the (fine) grid we 
collect the pressure.
Here we collect the pressure every four grid points, 
however exclude a 19 grid points deep layer around the boundary
(since the boundary conditions are known),
so that the number of measurement points is 49.
Collecting data this way allows us to use all
 data directly in our multiple 
grids approach with $16\times 16$ and $32\times 32$ grids 
(see below).
The data are perturbed with a Gaussian random variable
$r\sim\mathcal{N}(0,R)$, with a diagonal covariance matrix $R$
(i.e.~we assume that measurement errors are uncorrelated).
The variance at each measurement location is set
to 30\% of the reference solution.
This relatively large variance brings about significant 
non-Gaussian features in the posterior pdf.

\subsection{The log-normal prior, its discretization and dimensional reduction}
The prior for permeability fields is often assumed to be log-normal 
and we follow suit.  
The logarithm of the permeability $\kappa$ is thus a Gaussian field
with a squared exponential covariance function~\cite{Rasmussen2006},
\begin{equation}\label{equation:correlation}
R(x_1,x_2,y_1,y_2)=\sigma_x^2\sigma_y^2\exp \left(-\frac{\left(x_1-x_2\right)^2}{l_x^2}-\frac{\left(y_1-y_2\right)^2}{l_y^2}\right),
\end{equation}
where $(x_1,y_1)$, $(x_2,y_2)$ are two points in $\Omega$, 
and where the correlation length $l_x$ and $l_y$ 
and the parameters $\sigma_x,\sigma_y$ are given scalars. 
In the numerical experiments below, we choose $\sigma_x=\sigma_y=1$
and $l_x=l_y=\sqrt{0.5}$.
With this prior, 
we assume that the (log-) permeability is a smooth function of $x$ and $y$, 
so that solutions of~(\ref{equation:upeqn}) uniquely exist. 
Moreover, the theory presented in \cite{StuartInverse,Stuart2} applies 
and a well defined posterior also exists. 

We approximate the lognormal prior on the regular $N\times N$ grid 
by an $N^2$ dimensional log-normal random variable with
covariance matrix $\Sigma$ with elements $\Sigma(i,j)=R(x_i,x_j,y_i,y_j)$, 
$i,j=1,\dots,N$ where $N$ is the number of
grid points in each direction. 
To keep the computations manageable (for fine grids and large $N$),
we perform the matrix computations with a low-rank approximation of $\Sigma$
obtained via Karhunan-Lo\`{e}ve~(KL) expansions~\cite{Ghanem,LeMaitre}.
Specifically, the factorization of the covariance function $R(x_1,x_2,y_1,y_2)$ 
allows us to compute the covariance matrices in $x$ and $y$ directions separately, 
i.e.~we compute the matrices
\begin{equation*}
  \Sigma_x(i,j) = \sigma_x^2\exp \left(-\frac{(x_i-x_j)^2}{l_x^2}\right),\quad  \Sigma_y(i,j)=\sigma_y^2\exp \left(-\frac{(y_i-y_j)^2}{l_y^2}\right).
\end{equation*}
We then compute singular value decompositions (SVD) in each 
direction to form low-rank approximations $\hat{\Sigma}_x\approx \Sigma_x$ 
and $\hat{\Sigma}_y\approx \Sigma_y$ by neglecting small eigenvalues. 
These low rank approximations define a low rank approximation of the covariance matrix
\begin{equation*}
  \Sigma \approx \hat{\Sigma}_x \otimes \hat{\Sigma}_y,
\end{equation*}
where $\otimes$ is the Kronecker product.
Thus, the eigenvalues and eigenvectors of $\hat{\Sigma}$ are 
the products of the eigenvalues and eigenvectors of $\hat{\Sigma}_x$ and $\hat{\Sigma}_y$. 
The left panel of Figure 1
shows the spectrum of $\hat{\Sigma}$, 
and it is clear that the decay of the eigenvalues of $\Sigma$
is rapid and suggests a natural model reduction. 
\begin{figure} 
\label{fig:prior}
\begin{center}
\includegraphics[width=0.6\textwidth]{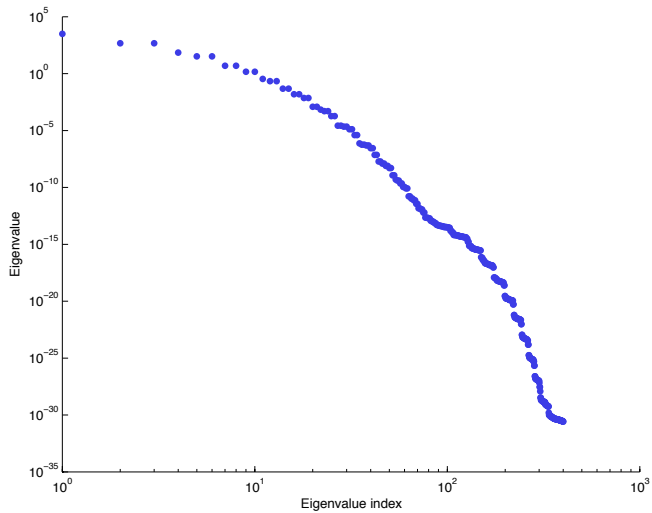}
\end{center}
\caption{Spectrum of the covariance matrix of lognormal prior.}
\end{figure}
If we neglect small eigenvalues (and set them to zero),
then
\begin{equation*}
  \hat{\Sigma} = V^T\Lambda V,
\end{equation*}
approximates $\Sigma$ (in a least squares sense in terms
of the Frobenius norms of $\Sigma$ and $\hat{\Sigma}$);
here $\Lambda$ is a diagonal matrix whose diagonal elements 
are the $m$ largest eigenvalues of $\Sigma$ and $V$ is an $m\times N$ 
matrix whose columns are the corresponding eigenvectors.
With $m=30$, we capture $99.9\%$ of the variance
(in the sense that the sum of the first 30 eigenvalues
is 99\% of the sum of all eigenvalues).
In reduced coordinates, the prior is
\begin{equation*}	
	\hat{K} \sim \ln \mathcal{N}\left(\hat{\mu},\hat{\Sigma}\right).
\end{equation*}
The linear change of variables
\begin{equation*}
	\theta = V^T\Lambda^{-0.5}\hat{K},
\end{equation*}
highlights that it is sufficient to estimate $m\ll N^2$ parameters 
(the remaining parameters are constrained by the prior).
We will carry out the computations in the reduced
coordinates $\theta$, for which the prior is 
\begin{equation}
\label{eq:PriorTheta}
  p(\theta)=\mathcal{N}\left(\mu,I_m\right),
\end{equation}
where $\mu= V^T\Lambda^{-0.5}\hat{\mu}$.
Note that this model reduction follows naturally 
from assuming that the
permeability is smooth,
so that errors correlate, 
and the probability mass localizes
in parameter space.
A similar observation, in connection with
data assimilation, was made in \cite{CM2013}.

Note that we achieve mesh-independence in the
implicit sampling algorithm we propose
by sampling in the $\theta$-coordinates
rather than in the physical coordinate system.
We consider this scenario because it
allows us to compare our approach with MCMC
that also samples in $\theta$-coordinates
and, therefore, also is mesh-independent.

\subsection{BFGS optimization with adjoints and multiple grids}
Implicit sampling requires minimization of $F$ in \EQ{FDef} which, 
for this problem and in reduced coordinates, takes the form
\begin{equation*}
F(\theta)= \frac{1}{2} \theta^T\theta+\frac{1}{2}\left(z-MP(\theta)\right)^TR^{-1}\left(z-MP(\theta)\right).
\end{equation*}
We solve the optimization problem using BFGS
coupled to an adjoint code to compute the gradient of $F$ with respect to $\theta$
(see also, e.g.~\cite{Oliver2011,Hinze}).

The adjoint calculations are as follows.
The gradient of $F$ with respect to $\theta$ is
\begin{equation*}
\nabla_{\theta} F(\theta)= \theta +\left(\nabla_{\theta} P(\theta)\right)^T W,
\end{equation*}
where $W=-M^TR^{-1}(z-MP(\theta))$. We use the chain rule to derive
$\left(\nabla_{\theta} P(\theta)\right)^T W $ as follows:
\begin{equation*}
(\nabla_{\theta} P(\theta))^TW=\left(\nabla_{K}
  P(\theta)\frac{\partial K}{\partial \hat{K}}\frac{\partial
    \hat{K}}{\partial \theta} \right)^T W
=\left(\nabla_{K}
  P(\theta) e^{\hat{K}} V\Lambda^{0.5}\right)^TW
=\left(V\Lambda^{0.5}\right)^T\left(\nabla_{K}
  P(\theta)e^{\hat{K}} \right)^T W,
\end{equation*}
where $e^{\hat{K}}$ is a $N^2\times N^2$ diagonal matrix whose elements are the exponentials of the components of $\hat{K}$.
The gradient $\nabla_{K} P(\theta)$ can be obtained directly from our
finite element discretization. Let $P=P(\theta)$ and let $K_{l}$ be the $l$th component of $K$, and take the derivative with respect to $K_{l}$ of \EQ{DEQ} to obtain
 \begin{equation*}
 \frac{\partial P}{\partial K_{l}  } =- A^{-1} 
 \frac{\partial A}{\partial K_{l}} P
\end{equation*} 
where $\partial A/\partial K_{l}$ are component-wise derivatives. We
use this result to obtain the following expression 
\begin{equation}\label{gradF}
\left(\nabla_{K}
  P(\theta) e^{\hat{K}}\right)^T W=-\left(e^{\hat{K}}\right)^T\left[ \begin{array}{c}
 P^T\frac{\partial A}{\partial K_{1}}\left(A^{-T}  W \right)\\
\vdots \\
   P^T\frac{\partial A}{\partial K_{N^2}} \left(A^{-T}  W \right)   
 \end{array} 
 \right].
\end{equation}
When $P$ is available, the most expensive part in (\ref{gradF}) is to evaluate $A^{-T}W$,
which is equivalent to solving the adjoint
problem (which is equal to itself for this self-adjoint problem). The
rest can be computed element-wise by the definition of $A$. Note that
there are only a  fixed number of nonzeros in each $\frac{\partial
  A}{\partial K_l}$, so that the additional work for solving the adjoint
problem in (\ref{gradF}) is about $O(N^2)$,
which is small compared to the work required for the adjoint solve.

 Collecting terms we finally obtain the gradient
\begin{eqnarray*}
\nabla_\theta F(\theta)= \theta+\left(V\Lambda^{0.5}\right)^T\left(\nabla_{K} P(\theta) e^{\hat{K}}\right)^TW 
&=&\theta-\left(V\Lambda^{0.5}\right)^T\left(e^{\hat{K}}\right)^T\left[ \begin{array}{c}
 P^T\frac{\partial A}{\partial K_{1}}\left(A^{-T}  W \right)\\
\vdots \\
   P^T\frac{\partial A}{\partial K_{N^2}} \left(A^{-T}  W \right)   
 \end{array} 
 \right]
.
\end{eqnarray*}
Multiplying by
$\left(V\Lambda^{0.5}\right)^T$ to
go back to physical coordinates will require an additional work of $O(mN^2)$. 
Note that the adjoint calculations for the gradient 
require only one adjoint solve because the forward solve 
(required for $P$) has already been done before the gradient calculation in the BFGS algorithm.
This concludes our derivation of an adjoint method for gradient computations.

The gradient is used in a BFGS method
with a cubic interpolation line search (see \cite[Chapter~3]{Nocedal}). 
We chose this method here because it defaults to taking the full step (of length 1) 
without requiring additional computations,
if the full step length satisfies the Wolfe conditions.
To reduce the number of fine-grid solves we use
the multiple grid approach described above with $16\times 16$, 
$32\times 32$ and $64\times 64$ grids.
We initialize the minimization on the course grid with the mean of the prior,
and observe a convergence after about 9 iterations, 
requiring 16 function and 16 gradient evaluations,
which corresponds to a cost of 32 coarse grid solves
(estimating the cost of adjoint solves with the cost of forward solves).
The result is used to initialize an optimization on a finer $32\times 32$  grid.
The optimization on $32\times 32$ grid converges in 6 iterations,
requiring 7 function and 7 gradient evaluations
(at a cost of 14 medium grid solves).
The solution on the medium grid is then used to initialize
 the finest $64\times 64$  grid optimization.
This optimization converges in 5 iterations, requiring
12 fine-grid solves.
We find the same minimum without the multiple grid approach,
i.e.~if we solve the minimization on the fine grid, 
however these computations require 36 fine grid solves.
The multiple-grids approach we propose
requires about 17 fine grid solves 
(converting the cost of coarse-grid solves to fine-grid solves)
and, thus, could significantly reduce
the number of required fine-grid solves.

\subsection{Implementation of  the random and linear maps}
Once the minimization is completed,
we generate samples using either the linear map
or random map methods described above.
Both require the Hessian of $F$ at
the minimum. 
We have not coded second-order adjoints,
and computing the Hessian with finite differences 
requires $m(m+1) = 930$ forward solutions,
which is expensive (and if $m$ becomes large,
this becomes infeasible).
Instead, we approximate the Hessian.
We found that the approximate Hessian of our BFGS
is not accurate enough to lead to a good implicit sampling method.
However, the Hessian approximation proposed in
\cite{Iglesias2} and often used in LMAP,
\begin{equation}
\label{eq:approxHess}
	H \approx \hat{H} = I - Q^T(QQ^T+R)^{-1}Q,
\end{equation}
where $Q = M\,\nabla_\theta  P$,
leads to good results (see below).
Here the gradient of the pressure (or the Jacobian)
is computed with finite differences,
which requires $m+1$ forward solves.
Note that the approximation is exact for linear Gaussian problems.

With the approximate Hessian, we define $L$
in (\ref{eq:LM}) and (\ref{equation:Ansatz})
as a Cholesky factor of $\hat{H}=LL^T$.
Generating samples with the random map method requires solving 
(\ref{equation:SamplingEquation}) with the ansatz (\ref{equation:Ansatz}). 
We use a Newton method for solving these equations
and observe that it usually converges quickly 
(within 1-4 iterations). 
Each iteration requires a derivative of $F$
with respect to $\lambda$, 
which we implement using the adjoint method,
so that each iteration requires two forward solutions.
In summary, the random map method requires 
between 2-8 forward solutions per sample.
 
The linear map method requires generating a sample
using (\ref{eq:LM}) and weighing it by (\ref{eq:weightsLinear}).
Evaluation of the weights thus requires one forward solve.
Neglecting the cost for the remaining linear algebra,
the linear map has a cost of 1 PDE solve per sample.
 
We assess the quality of the weighted samples 
by the variance of the weights: 
the sampling method is good if the 
variance of the weights is small.
In particular, if the weights are constant, 
then this variance is zero and the sampling method is perfect.
The variance of the weights is equal to $R-1$, where
\begin{equation*}
	R = \frac{E(w^2)}{E(w)^2}.
\end{equation*}
In fact, $R$ itself can be used to measure the quality
of the samples \cite{VEW12,AMGC02}.
If the variance of the weights is small,
then $R\approx 1$.
Moreover, the effective sample size,
i.e.~the number of unweighted samples that would be 
equivalent in terms of statistical accuracy to the
set of weighted samples,
is about $M/R$ \cite{VEW12},
where $M$ is the number of samples we draw.
In summary, an $R$ close to one indicates 
a well distributed weighted ensemble.

We evaluate $R$ for 10 runs with $M=10^4$ samples
for each method and find, for the linear map method,
a mean of $1.79$ and standard deviation $0.014$,
and for the random map method a mean of $1.77$
and standard deviation $0.013$.
The random map method thus performs slightly better,
however the cost per sample is also slightly larger
(because generating a sample requires solving 
(\ref{equation:SamplingEquation}), which in turn
requires solving the forward problem).
Because the linear map method produces weighted ensembles of
about the same quality as the random map,
and since the linear map is less expensive and easier to program,
we conclude that the linear map is a more natural choice for this example.
 
We have also experimented with symmetrization of implicit sampling~\cite{GLM}.
The idea is similar to the classic trick of antithetic variates~\cite{Kalos}.
The symmetrization of the linear map method is as follows.
Sample the reference variable to obtain a $\xi$ and compute a sample $x^+$ using~(\ref{eq:LM}).
Use the same $\xi$ to compute $x^- = \mu-L^{-T}\xi$.
Then pick $x$ with probability $p^+ = w(x)/(w(x)+w(x^-))$
and pick $x^-$ with probability $p^- =w(x^-)/(w(x)+w(x^-))$,
and assign the weight $w^s = (w(x^+)+w(x^-))/2$.
This symmetrization can lead to a smaller $R$,
i.e.~a better distributed ensemble, in the small noise limit.
In our example, we compute the quality measure $R$ of 1.67.
While this $R$ is smaller than for the non-symmetrized methods,
the symmetrization does not pay off in this example,
since each sample of the symmetrized method requires two
forward solves (to evaluate the weights).
 
Note that we neglect computations other than the forward model
evaluations when we estimate the computational cost of the sampling algorithms
(as we did with the BFGS optimization as well).
This is justified because computations with $\theta$
(e.g.~generating a sample using the linear map method)
is inexpensive due to the model reduction via Karhunen-Lo\`{e}ve. 
 
We illustrate how our algorithms can be used 
by presenting results of a typical run and for a typical problem set-up
in terms of e.g.~strength of the observation noise and the number of observations.
We tested our algorithms in a variety of other settings as well,
and observed that our methods operate reliably in different problem set-ups,
however found that many of the problems one can set up are 
almost Gaussian problems and therefore easy to solve.
We present here a case where the large observation noise (see above)
brings about significant non-Gaussian features in the posterior.

Shown in Figure 2 are the true permeability 
(the one we use to generate the data) on the left,
the mean of the prior in the center,
and the conditional mean we computed with the linear map
method and $10^4$ samples on the right.
\begin{figure} 
\label{fig:result}
\begin{center}
\includegraphics[width=1\textwidth]{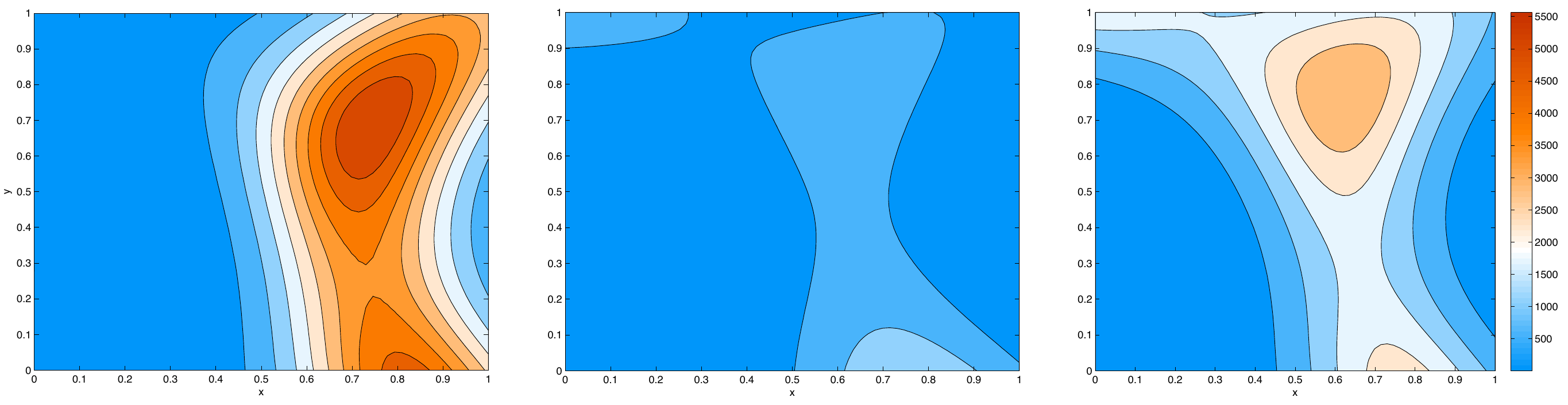}
\end{center}
\caption{Left: true permeability that generated the data.
Center: mean of prior. 
Right: conditional mean computed with implicit sampling with random maps.}
\end{figure}
We observe that the prior is not very informative,
in the sense that it underestimates the permeability considerably.
The conditional mean captures most of the large scale features,
such as the increased permeability around $x=0.7$, $y=0.8$,
however, there is considerable uncertainty in the posterior.
Figure 3 illustrates this uncertainty
and shows four samples of the posterior.
\begin{figure} 
\begin{center}
\includegraphics[width=0.8\textwidth]{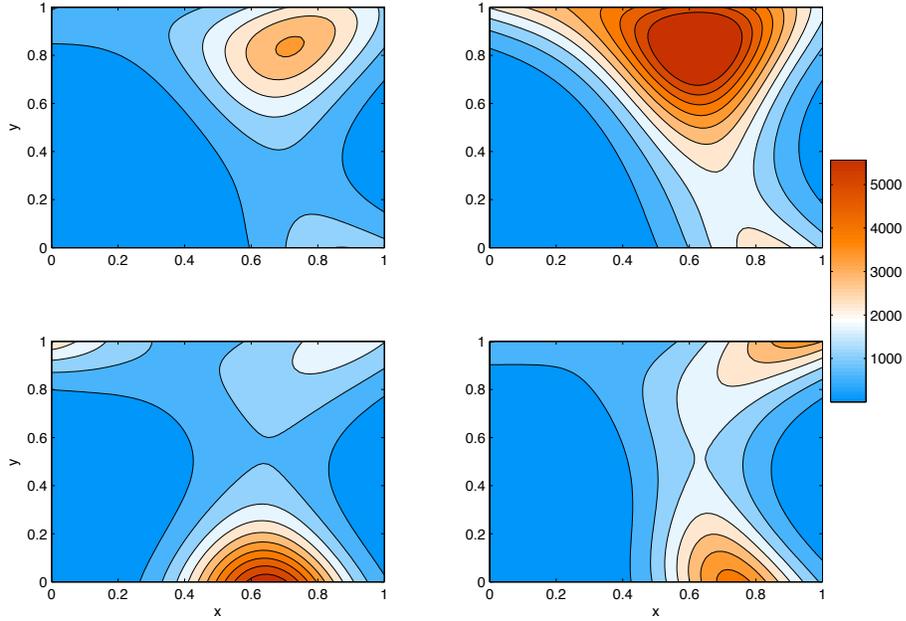}
\end{center}
\label{fig:PosteriorSamples}
\caption{Four samples of the posterior generated by implicit sampling with random maps.}
\end{figure}
These samples are obtained by resampling the weighted ensemble,
so that one is left with an equivalent unweighted set of samples,
four of which are shown in Figure 4.
The four samples correspond to rather different subsurface structures. 
If more accurate and more reliable estimates of the permeability
are required, one must increase the resolution of the data
or reduce the noise in the data.
 
\subsection{Connections with other methods}
We discuss connections of our implicit sampling schemes
with other methods that are in use in subsurface flow 
parameter estimation problems. 

\subsubsection{Connections with linearization about the MAP}
One can estimate parameters by computing the MAP point,
i.e.~the most likely parameters in view of the data
\cite{Iglesias2,Oliver2011}.
This method, sometimes called the MAP method,
can make use of the multiple grids approach presented here, 
however represents an incomplete solution to the Bayesian parameter estimation problem, 
because the uncertainty in the parameters may be large,
however the MAP point itself contains no information
about this uncertainty.
To estimate the uncertainty of the MAP point,
one can use linearization about the MAP point (LMAP)
\cite{Iglesias2,Oliver2011,Dean2007,Bui},
in which one computes the MAP point and the Hessian
of $F$ at the MAP point and uses the inverse of this Hessian
as a covariance.
The cost of this method is the cost
of the optimization plus the cost of computing the Hessian.
For the example above,
LMAP overestimates the uncertainty
and gives a standard deviation of $0.61$
for the first parameter $\theta_1$.
The standard deviation we compute with the
linear map and random map methods however is $0.31$.
The reason for the over-estimation of the uncertainty with LMAP
is that the posterior is not Gaussian.
This effect is illustrated in Figure~4
where we show histograms of the marginals
of the posterior for the first four parameters
$\theta_1,\theta_2,\theta_3,\theta_4$,
along with their Gaussian approximation as in LMAP.
\begin{figure} 
\begin{center}
\includegraphics[width=0.6\textwidth]{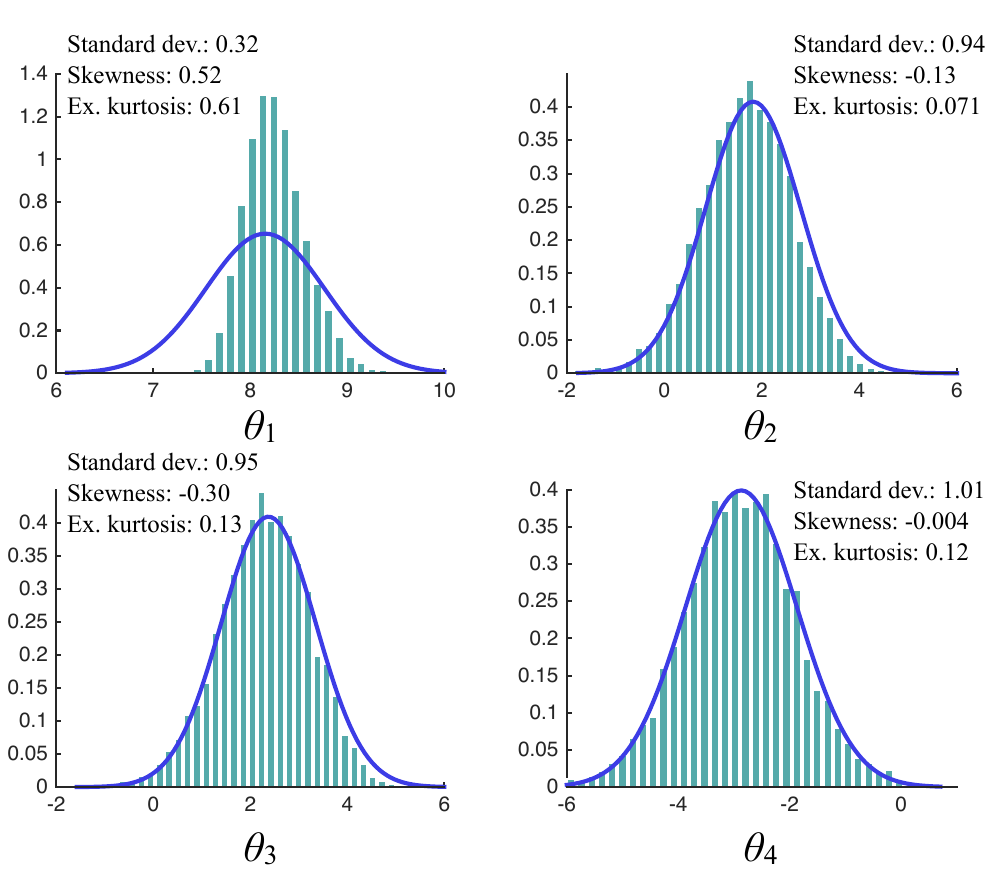}
\end{center}
\label{fig:Gaussian}
\caption{Marginals of the posterior computed with
implicit sampling with random maps and their Gaussian
approximation obtained via LMAP.
Top left: $p(\theta_1|z)$.
Top right: $p(\theta_2|z)$.
Bottom left: $p(\theta_3|z)$.
Bottom right: $p(\theta_4|z)$.
}
\end{figure}
We also compute the skewness and excess kurtosis
for these marginal densities.
While the marginals for the parameters may become
``more Gaussian'' for the higher order coefficients
of the KL expansion, the joint posterior exhibits significant non-Gaussian behavior.
Since implicit sampling (with random or linear maps)
does not require linearizations
or Gaussian assumptions, it can correctly capture these
non-Gaussian features.
In the present example, accounting for the non-Gaussian\slash
nonlinear effects brings about a reduction of the uncertainty
(as measured by the standard deviation) by a factor of two in the parameter $\theta_1$.

Note that code for LMAP, 
can be easily converted into an implicit sampling code.
In particular, implicit sampling with linear maps
requires the MAP point and an approximation of the Hessian
at the minimum. Both can be computed with LMAP codes.
Non-Gaussian features of the posterior can then be captured
by weighted sampling, where each sample comes at a cost
of a single forward simulation.

\subsubsection{Connections with Markov Chain Monte Carlo}
Another important class of methods
for solving Bayesian parameter estimation problems is MCMC 
(see e.g.~\cite{Iglesias2} for a discussion of MCMC in subsurface flow problems).
First we compare implicit sampling with Metropolis MCMC \cite{Liu2008},
where we use an isotropic Gaussian proposal density,
for which we tuned the variance to achieve an acceptance rate of about 30\%.
This method requires one forward solution per step
(to compute the acceptance probability).
We start the chain at the MAP (to reduce burn-in time).
In figure~5 we show the approximation of the conditional mean
of the variables $\theta_1,\theta_2$, and $\theta_5$,
as a function of the number of steps in the chain (left)
to illustrate the behavior of the MCMC chain.
\begin{figure} 
\begin{center}
\includegraphics[width=1\textwidth]{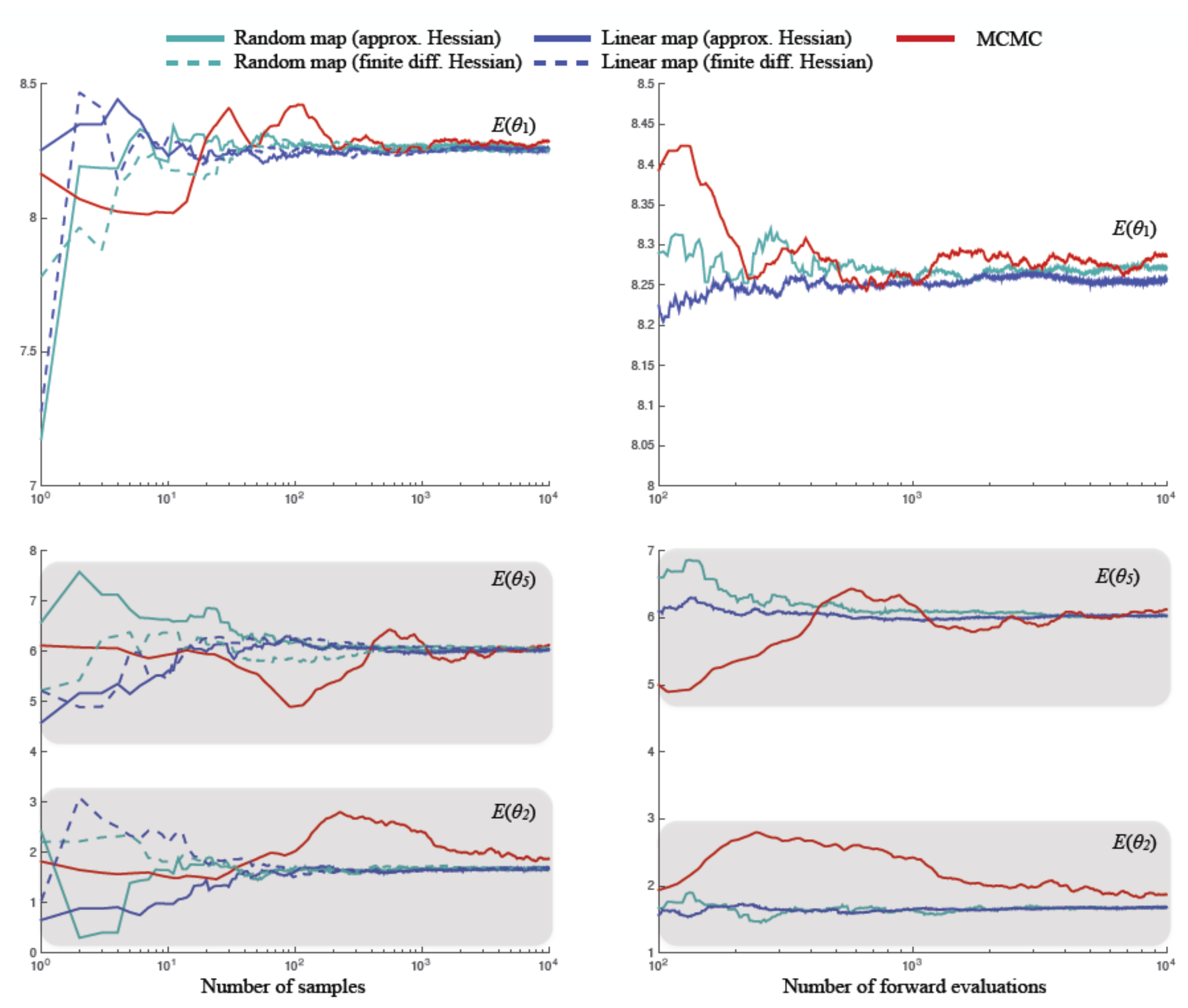}
\end{center}
\caption{Expected value as a function of the number of samples (left),
and as a function of required forward solves (right).
Red: MCMC. Turquoise: implicit sampling with random maps
and approximate Hessian (dashed) and finite difference Hessian (solid).
Blue: implicit sampling with linear maps
and approximate Hessian (dashed) and finite difference Hessian (solid).}
\end{figure}
We observe that, even after $10^4$ steps,
the chain has not settled,
in particular for the parameter $\theta_2$ (see bottom left).
With implicit sampling we observe a faster convergence,
in the sense that the approximated conditional mean does not
change significantly with the number of samples.
In fact, about $10^2$ samples are sufficient for accurate estimates
of the conditional mean.
As a reference solution, we also show results we obtained
with implicit sampling (with both random and linear maps)
for which we used a Hessian computed with
finite differences (rather than with the approximation in
equation~(\ref{eq:approxHess})).

The cost per sample of implicit sampling 
and the cost per step of Metropolis MCMC are different,
and a fair comparison of these methods should take
these costs into account.
In particular, the off-set cost of the minimization and
computation of the Hessian, required for implicit sampling 
must be accounted for.
We measure the cost of the algorithms by the number of
forward solves required (because all other computations
are negligible due to the model reduction).
The results are shown for the parameters $\theta_1,\theta_2$
and $\theta_5$ in the right panels of Figure~5.

We find that the fast convergence 
of implicit sampling makes up for the
relatively large a priori cost
(for minimization and Hessian computations).
In fact, the figure suggests that the random and linear map methods require
about $10^3$ forward solves while Metropolis MCMC converges slower
and shows significant errors even after running the chain for $10^4$ steps,
requiring $10^4$ forward solves.
The convergence of Metropolis MCMC can perhaps be increased
by further tuning, or by choosing a more advanced transition density.
Implicit sampling on the other hand requires little tuning
other than deciding on standard tolerances for the optimization.
Moreover, implicit sampling generates independent samples
with a known distribution, so that issues such as
determining burn-in times, auto-correlation times and 
acceptance ratios, do not arise.
It should also be mentioned that implicit sampling
is easy to parallelize (it is embarrassingly parallel
once the optimization is done). 
Parallelizing Metropolis MCMC on the other hand is not trivial,
because it is a sequential technique.

Finally, we discuss connections of our proposed implicit sampling methods to a new MCMC 
method, stochastic Newton MCMC \cite{Martin}.  
In stochastic Newton one first finds the MAP point (as in implicit sampling)
and then starts a number of MCMC chains from the MAP point.
The transition probabilities are based on local information about 
$F$ and make use of the Hessian of $F$, evaluated at the location of the chain.
Thus, at each step, a Hessian computation is required which,
with our finite difference scheme, 
requires 31 forward solves (see above) and, therefore,
is expensive (compared to generating samples with implicit sampling,
which requires computing the Hessian only once).
Second-order adjoints (if they were available)
do not reduce that cost significantly.
We have experimented with stochastic Newton in our example
and have used 10--50 chains
and taking about 200 steps per chain.
Without significant tuning,
we find acceptance rates of only a few percent,
leading to a slow convergence of the method.
We also observe that the Hessian may not be positive definite
at all locations of the chain and, therefore, 
can not be used for a local Gaussian transition probability. 
In these cases, we use a modified Cholesky algorithm
(for affine invariance) to obtain a definite matrix that can be used 
as a covariance of a Gaussian.
In summary, we find that stochastic Newton MCMC
is impractical for this example because
the cost of computing the Hessian is too large with our
finite differences approach.
Variations of stochastic Newton were explained in~\cite{Petra}.
The stochastic Newton MCMC with MAP-based Hessian
is the stochastic Newton method as above,
however the Hessian is computed only at the MAP point
and then kept constant throughout the chain.
The ``independence sampling with a MAP point-based Gaussian proposal'' (ISMAP)
is essentially an MCMC version of the linear map method described above.
The ISMAP MCMC method is to use the Gaussian approximation
of the posterior probability at the MAP point as the proposal density for MCMC.
The samples are accepted or rejected based on the weights of the linear map method
described above.
ISMAP is also easier to parallelize than stochastic Newton or stochastic Newton
with MAP-based Hessian.

\section{Conclusions}
We explained how to use implicit sampling 
to estimate the parameters in PDE 
from sparse and noisy data. 
The idea in implicit sampling is to find
the most likely state, often called
the maximum a posteriori (MAP) point,
and generate samples that explore the 
neighborhood of the MAP point.
This strategy can work well if
the posterior probability mass localizes around the MAP point,
which is often the case when the data constrain the parameters.
We discussed how to implement these ideas efficiently
in the context of parameter estimation problems.
Specifically, we demonstrated that our approach is 
mesh-independent in the sense that we sample
finite dimensional subspaces even when the grid is refined.
We further showed how to use multiple grids
to speed up the required optimization,
and how to use adjoints for the optimization and during sampling.

Our implicit sampling approach has the advantage that
it generates independent samples, so that issues connected with MCMC,
e.g.~estimation of burn-in times, auto-correlations of the samples, or tuning of
acceptance ratios, are avoided. 
Our approach is also fully nonlinear and captures non-Gaussian
features of the posterior (unlike linear methods such as
the linearization about the MAP point) and is easy to parallelize.

We illustrated the efficiency of our approach
in numerical experiments with an elliptic inverse problem
that is of importance in applications to 
reservoir simulation\slash management and pollution modeling. 
The elliptic forward model is discretized using finite elements,
and the linear equations are solved by balancing domain decomposition by constraints. 
The optimization required by implicit sampling is done with
with a BFGS method coupled to an adjoint code for gradient calculations.
We use the fact that the solutions are expected to be smooth 
for model order reduction based on Karhunan-Lo\`{e}ve expansions, 
and found that our implicit sampling approach can exploit this low-dimensional structure.
Moreover, implicit sampling is about an order of magnitude
faster than Metropolis MCMC sampling (in the example we consider). 
We also discussed connections and differences of our approach
with linear\slash Gaussian methods, such as linearization about the MAP,
and with stochastic Newton MCMC methods.

\section*{Acknowledgements}
This material is based upon work supported by the 
U.S.~Department of Energy, Office of Science,
Office of Advanced Scientific Computing Research, 
Applied Mathematics program under contract DE-AC02005CH11231, 
and by the National Science Foundation under grant DMS-0955078, DMS-1115759,
DMS-1217065, and DMS-1419069.

\bibliographystyle{plain}

\begin{thebibliography}{10}

\bibitem{AMGC02}
M.S.~Arulampalam, S.~Maskell, N.~Gordon, and T.~Clapp.
\newblock A tutorial on particle filters for online nonlinear\slash non-Gaussian
Bayesian tracking
\newblock {\em IEEE Transactions on Signal Processing}, 50(2):174--188, 2002.

\bibitem{amc}
E.~Atkins, M.~Morzfeld, and A.J.~Chorin.
\newblock Implicit particle methods and their connection with variational data assimilation.
\newblock {\em Monthly Weather Review}, 141(6):1786--1803, 2013.




\bibitem{Bear}
J.~Bear.
\newblock {\em Modeling groundwater flow and pollution}.
\newblock Kluwer, 1990.

\bibitem{Blb}
P.~Bickel, B.~Li, and T.~Bengtsson.
\newblock Sharp failure rates for the bootstrap particle filter in high
  dimensions.
\newblock {\em IMS Collections: Pushing the Limits of Contemporary Statistics:
  Contributions in Honor of Jayanta K. Ghosh}, 3:318--329, 2008.
  
\bibitem{Blb2}
T.~Bengtsson, P.~Bickel, and B.~Li.
\newblock Curse of dimensionality revisited: the collapse of importance sampling in very large scale systems.
\newblock {\em MS Collections: Probability and Statistics: Essays in Honor of
David A. Freedman}, 2:316--334, 2008.


\bibitem{Bui}
T.~Bui-Thanh, O.~Ghattas, J,~Martin, and G.~Stadler.
\newblock A computational framework for infinite-dimensional {B}ayesian
              inverse problems {P}art {I}: {T}he linearized case, with
              application to global seismic inversion
\newblock {\em SIAM J. Sci. Comput.}, 35(6):A2494-A2523, 2013.

\bibitem{Braess}
D. Braess.
\newblock Finite Elements: Theory, Fast Solvers, and
                  Applications in Solid Mechanics.
  \newblock Cambridge University Press, 1997.

\bibitem{ChorinHald}
A.J.~Chorin and O.H.~Hald.
\newblock {\em Stochastic Tools in Mathematics and Science}.
\newblock Springer, third edition, 2013.

\bibitem{cmt}
A.J.~Chorin, M.~Morzfeld, and X.~Tu.
\newblock Implicit particle filters for data assimilation.
\newblock {\em Communications in Applied Mathematics Computational Sciences}, 5:221--240, 2010.

\bibitem{cmt2}
A.J.~Chorin, M.~Morzfeld, and X.~Tu.
\newblock Implicit sampling, with applications to data assimilation.
\newblock {\em Chinese Annals of Mathematics}, 34B:89--98, 2013.

\bibitem{CM2013}
A.J.~Chorin and M.~Morzfeld.
\newblock Condition for successful data assimilation.
\newblock {\em Journal of Geophysical Research}, 118(20):11522--11533, 2013.

\bibitem{CT1}
A.J.~Chorin and X.~Tu.
\newblock Implicit sampling for particle filters.
\newblock {\em Proceedings of the National Academy Sciences USA}, 106:17249--17254, 2009.

\bibitem{Stuart2}
M.~Dashti and A.M.~Stuart.
\newblock Uncertainty quantification and weak approximation of an elliptic
  inverse problem.
\newblock {\em SIAM Journal on Numerical Analysis}, 49(6):2524--2542, 2011.

\bibitem{BDDC}
C.~R.~Dohrmann.
\newblock A preconditioner for substructuring based on constrained energy
          minimization.
\newblock {\em SIAM Journal on Scientific Computing}, 25:246--258, 2003.
 
\bibitem{Doucet2000}
A.~Doucet, S.~Godsill, and C.~Andrieu.
\newblock On sequential {M}onte {C}arlo sampling methods for {B}ayesian filtering.
\newblock {\em Statistics and Computing}, 10:197--208, 2000.

\bibitem{Efendiev}
Y.~Efendiev, T.~Hou, and W.~Luo.
\newblock Preconditioning {M}arkov chain {M}onte {C}arlo simulations using
  coarse-scale models.
\newblock {\em SIAM Journal on Scientific Computing}, 28(2):776--803 (electronic), 2006.

\bibitem{Fedorenko}
R.P.~Fedorenko.
\newblock A relaxation method for solving elliptic difference equations.
\newblock {\em USSR Computational Mathematics and Mathematical Physics},1, 1961.


\bibitem{Ghanem}
R.~Ghanem and P.~Spanos.
\newblock {\em Stochastic Finite Elements: A Spectral Approach}.
\newblock Dover, 2003.

\bibitem{Goodman}
J.~Goodman and A.D.~Sokal.
\newblock Multigrid Monte Carlo method. Conceptual foundations.
\newblock {\em Physical Review D}, 40:2035--2071, 1989.

\bibitem{GLM}
J.~Goodman and K.L.~Lin and M.~Morzfeld.
\newblock Small-noise analysis and symmetrization of implicit Monte Carlo samplers.
\newblock {\em Communications on Pure and Applied Mathematics}, accepted, 2015.


\bibitem{Hinze}
M.~Hinze, R.~Pinnau, M.~Ulrbich, and S.~Ulbrich.
\newblock {\em Optimization with PDE Constraints}.
\newblock Springer, 2009.




\bibitem{Iglesias}
M.A.~Iglesias, K.J.H.~Law and A.M.~Stuart.
\newblock Ensemble Kalman methods for inverse problems.
\newblock {\em Inverse Problems}, 29:045001(20 pp.), 2013.

\bibitem{Iglesias2}
M.A.~Iglesias, K.J.H.~Law and A.M.~Stuart.
\newblock Evaluation of {G}aussian approximations for data
assimilation in reservoir models.
\newblock {\em Computational Geosciences},  17: 851-885, 2013.

\bibitem{Kalos}
M.~Kalos and P.~Whitlock. 
\newblock {\em {M}onte {C}arlo methods, volume 1}.
\newblock John Wiley \& Sons, 1 edition, 1986.

\bibitem{LeMaitre}
O.P.~LeMaitre and O.M.~Knio. 
\newblock {\em Spectral Methods for Uncertainty Quantification: with Applications to Computational Fluid Dynamics}.
\newblock Springer, 2010.

\bibitem{Liu2008}
J.S.~Liu. 
\newblock {\em Monte Carlo Strategies for Scientific Computing}.
\newblock Springer, 2008.

\bibitem{liuchen1995}
J.S.~Liu and R.~Chen.
\newblock Blind Deconvolution via Sequential Imputations.
\newblock {\em Journal of the American Statistical Association}, 90(430):567--576, 1995.

\bibitem{Martin}
J.~Martin, L.C.~Wilcox, C.~Burstedde, and O.~Ghattas.
\newblock A stochastic {N}ewton {MCMC} method for large-scale
              statistical inverse problems with application to seismic
              inversion.
\newblock{\em SIAM J. Sci. Comput.} 34:A1460--A1487, 2012.


\bibitem{Marzouk1}
Y.M.~Marzouk, H.N.~Najm, and L.A.~Rahn.
\newblock Stochastic spectral methods for efficient Bayesian solution of inverse problems.
\newblock {\em Journal of Computational Physics}, 224(2):560--586, 2007.

\bibitem{Morzfeld2012}
M.~Morzfeld and A.J.~Chorin.
\newblock Implicit particle filtering for models with partial noise, and an
  application to geomagnetic data assimilation.
\newblock {\em Nonlinear Processes in Geophysics}, 19:365--382, 2012.

\bibitem{Morzfeld2011}
M.~Morzfeld, X.~Tu, E.~Atkins, and A.J.~Chorin.
\newblock A random map implementation of implicit filters.
\newblock {\em Journal of Computational Physics}, 231(4):2049--2066, 2012.

\bibitem{MarzoukOptimalMaps}
T.A.~Moselhy and Y.M.~Marzouk.
\newblock Bayesian inference with optimal maps.
\newblock {\em Journal of Computational Physics}, 231:7815--7850, 2012.

\bibitem{Nocedal}
J.~Nocedal and S.T.~Wright.
\newblock {\em Numerical Optimization}.
\newblock Springer, second edition, 2006.

\bibitem{Oliver14}
D.S.~Oliver.
\newblock Minimization for conditional simulation.
\newblock {\em Journal of Computational Physics}, 265:1--15, 2014. 

\bibitem{Dean2007}
D.S.~Oliver, A.C.~Reynolds, and N.~Liu.
\newblock {\em Inverse theory for petroleum reservoir characterization and history matching}.
\newblock Cambridge University Press, 2008.

\bibitem{Oliver2011}
D.S.~Oliver and Y.~Chen.
\newblock Recent progress on reservoir history matching: a review.
\newblock {\em Computers and Geosciences}, 15, 185--221, 2011.

\bibitem{Petra}
N.~Petra, J.~Martin, G.~Stadler, and O.~Ghattas,.
\newblock A computational framework for infinite-dimensional {B}ayesian
              inverse problems {P}art {II}: {S}stochastic Newton MCMC with
              application to ice sheet flow inverse problems
\newblock {\em SIAM J. Sci. Comput.}, 36(4):A1525-A1555, 2014.

\bibitem{Rasmussen2006}
C.E.~Rasmussen and C.K.I.~Wiliams.
\newblock {\em Gaussian processes for machine learning}.
\newblock MIT Press, 2006.

\bibitem{Tabak}
N.~Recca.
\newblock A new methodology for importance sampling.
\newblock{\em Masters Thesis}, 
Courant Institute of Mathematical Sciences, New York University.


\bibitem{Sny}
C.~Snyder, T.~Bengtsson, P.~Bickel, and J.~Anderson.
\newblock Obstacles to high-dimensional particle filtering.
\newblock {\em Monthly Weather Review}, 136:4629--4640, 2008.

\bibitem{StuartInverse}
A.~M. Stuart.
\newblock Inverse problems: a Bayesian perspective.
\newblock {\em Acta Numerica}, 19:451--559, 2010.




\bibitem{VEW12}
E.~Vanden-Eijnden and J.~Weare.
\newblock Data assimilation in the low-noise, accurate observation regime
with application to the Kuroshio current.
\newblock {\em Monthly Weather Review}, 141:1822--1841, 2012.

\bibitem{Brad}
B.~Weir, R.N.~Miller and Y.H.~Spitz.
\newblock Implicit estimation of ecological model parameters.
\newblock {\em Bulletin of Mathematical Biology}, 75:223--257, 2013.

\bibitem{Zaritski}
V.S.~Zaritskii and L.I.~Shimelevich.
\newblock Monte Carlo technique in problems of optimal data processing.
\newblock {\em Automation and Remote Control}, 12:95--103, 1975.


\end{thebibliography}

\end{document}